\documentclass[review]{elsarticle}

\usepackage{lineno,hyperref,amsmath,amsfonts,bm,algorithmic,algorithm,threeparttable,graphicx,subfig,cleveref,color}
\modulolinenumbers[5]
\usepackage{amsopn}

\DeclareMathOperator{\argmax}{argmax}
\crefname{algorithm}{Algorithm}{Algorithms}
\crefname{equation}{Eq.}{Eq.}
\Crefname{equation}{Equation}{Equations}
\crefname{table}{Table}{Tables}
\crefname{figure}{Fig.}{Figures}
\Crefname{figure}{Figure}{Figures}
\newcommand*{\reviewerA}[1]{\textcolor{black}{#1}}
\journal{Journal of Mechanical Systems and Signal Processing}

\begin{document}

\begin{frontmatter}

\title{Fast Greedy Optimization of Sensor Selection in Measurement with Correlated Noise
\tnoteref{mytitlenote}}
\tnotetext[mytitlenote]{Full document are available as \href{https://doi.org/10.1016/j.ymssp.2021.107619}{10.1016/j.ymssp.2021.107619}.}

\author[tohoku]{Keigo Yamada\corref{corresponding}
}
\cortext[corresponding]{Corresponding author}
\ead{yamada.keigo@aero.mech.tohoku.ac.jp}
\author[tohoku]{Yuji Saito}
\author[tohoku]{Koki Nankai}
\author[tohoku]{Taku Nonomura}
\author[tohoku]{Keisuke Asai}
\author[nagoya]{Daisuke Tsubakino}


\address[tohoku]{Department of Aerospace Engineering, Tohoku University, Sendai, Japan}
\address[nagoya]{Department of Aerospace Engineering, Nagoya University, Nagoya, Japan}

\begin{abstract}
A greedy algorithm is proposed for sparse-sensor selection in reduced-order sensing that contains correlated noise in measurement. 
The sensor selection is carried out by maximizing the determinant of the Fisher information matrix in a Bayesian estimation operator.
The Bayesian estimation with a covariance matrix of the measurement noise and a prior probability distribution of estimating parameters, which are given by the modal decomposition of high dimensional data, robustly works even in the presence of the correlated noise.  
After computational efficiency of the algorithm is improved by a low-rank approximation of the noise covariance matrix, the proposed algorithms are applied to various problems.
The proposed method yields more accurate reconstruction than the previously presented method with the determinant-based greedy algorithm, with reasonable increase in computational time.
\end{abstract}

\begin{keyword}
Data processing, sensor placement optimization, greedy algorithm, Bayesian state estimation.
\MSC[2020]     
    \href{https://zbmath.org/classification/?q=cc\%3A37M05}{37M05}\sep
   \href{https://zbmath.org/classification/?q=cc\%3A90C27}{90C27}\sep
   \href{https://zbmath.org/classification/?q=cc\%3A62F15}{62F15}\sep
\end{keyword}

\end{frontmatter}

\section{Introduction}

\reviewerA{Monitoring complex fluid behavior is essential for effective feedback control in aerospace engineering~\cite{kincaid2002d, corke2007sdbd}.
However, a number of sensors deployed in the field is limited due to a requirement for processing in the real-time situation and reducing communication energy. 
Therefore, optimization methods for sensor location and practical state-estimation schemes are critical, since one must get over a difficulty in estimation of off-wing fast fluid flows from a handful of sensors on wings.
Similar configurations assuming sparsity in acquired data and reducing amount of sensors are found in optimal experiment design~\cite{bates1996experimental}, compressed sensing~\cite{brunton2014compressive} and system identification~\cite{udwadia1994methodology, brunton2016discovering}. }

\reviewerA{Approaches for sensor selection varied very widely as seen in recent monographs, {\it e.g.} optimization after equation-based modeling~\cite{chaturantabut2010nonlinear} and machine learned optimization that circumvents modeling~\cite{semaan2017optimal}. However, this work presents optimization on sensor placement through the modeling of a phenomenon, which is constructed in a data-driven manner seen in Ref.~\cite{berkooz1993proper}. The reason is that high-dimensional visualization data can be obtained in the experiments of fluid dynamics or others by the recent development of measurement systems, and that previously mentioned requirement for feedback control motivates sensor selection in that way.
There are various objective metrics and optimization methods with regard to sensor selection algorithms.  The objective functions for the sensor selection problem are exemplified by those associated with the Fisher information matrix in ordinary linear least squares estimation~\cite{peherstorfer2018stability, nakai2020effect} or with a steady state error covariance matrix of the Kalman filter~\cite{ye2018complexity}. On the other hand, the optimization methods (as heuristic ones to the brute-force search) for the sensor selection problems are based on convex relaxation methods~\cite{joshi2009sensor,nonomura2020randomized}, greedy methods~\cite{manohar2018data,saito2019determinantbased} and proximal optimization algorithms~\cite{dhingra2014admm,nagata2020data-driven}.}

\reviewerA{In these established methodologies, however, formulations that considered spatially correlated noise could rarely be seen. The noise becomes often problematic, due to difference between a model and a phenomenon itself, in processing of acoustic signals~\cite{o2016distributed} and vibration~\cite{castro2013robustness}, and data-driven reduced-order modeling~\cite{saito2019determinantbased}.}
\reviewerA{Liu {\it et al.}~\cite{liu2016sensor} pointed out a simplified assumption of weakly correlated noise of Ref.~\cite{joshi2009sensor} and introduced formulations based on the trace of the Fisher information matrix with a general kernel of noise covariance between sensors. Uci\'{n}ski~\cite{ucinski2020d} developed an iterative optimization method for a similar objective function involving nonconvex terms to promote the Fisher information matrix to be regular in the relaxed form.}
%

Here, the aim of this paper is to improve sparse sensing and sensor selection algorithms \reviewerA{by several contributions provided in the list below;
\begin{itemize}
    \item Introduces a data-driven noise covariance matrix and Bayesian priors, that are both calculated from proper orthogonal decomposition (POD) on data. In our framework, modes generated by the POD procedure are divided into two; first $r$ modes corresponding to a order-reduced phenomenon and the rest corresponding to the correlated noise.
    \item Presents an objective function based on the determinant of the Fisher information matrix of a Bayesian state estimation and a greedy algorithm that leverages rank-one lemma as proposed in Ref.~\cite{saito2019determinantbased}.
    \item Develops an efficient implementation by involving an approximation of the noise covariance matrix and verifies its effect on numerical simulations and several actual datasets.
\end{itemize}
The approach and the contribution of the present study are illustrated in \cref{fig:flowchart}.
\begin{figure}[htbp]
    \centering
    \includegraphics[width=4in]{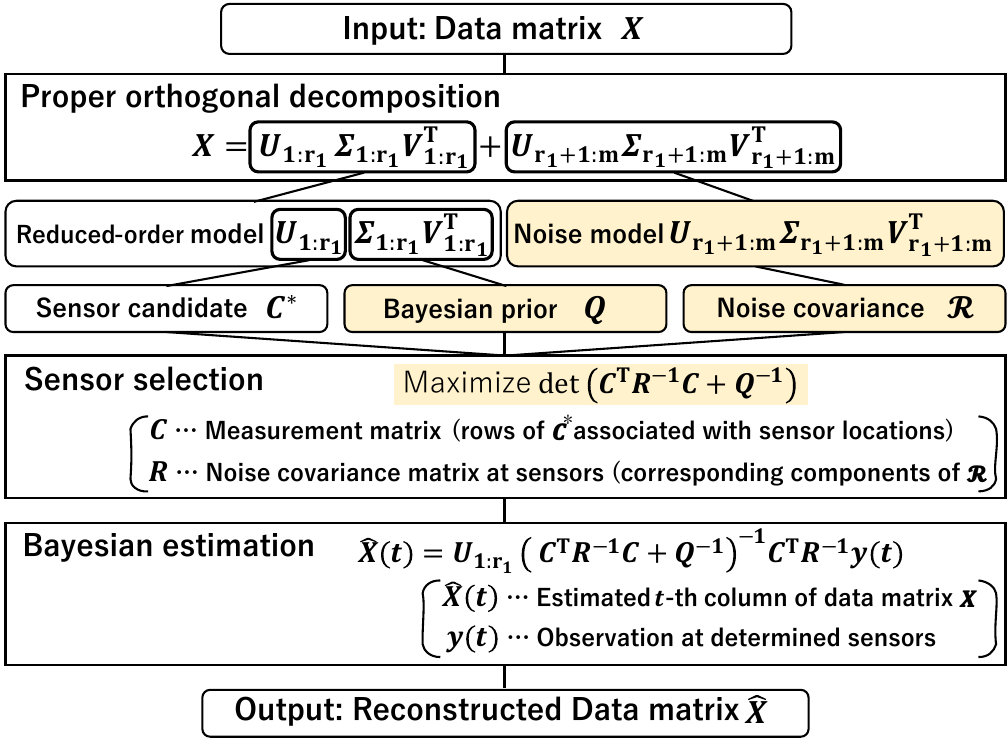}
    \caption{Concept of the manuscript (based on \cref{alg_BDG}). Highlighted steps are introduced in the present study. This flowchart reduces to that of the previous framework by ignoring highlighted steps and by replacing ``Bayesian estimation'' by ``least-squares estimation.''}\label{fig:flowchart}
\end{figure}
}
Firstly, basics of the POD-based reduced-order modeling and sparse sensing are briefly revisited. Algorithms of the previous and the present study are given in \cref{subsection:basics,subsection:proposed}, and then the superiority of the proposed method for noisy datasets is shown by reproducing randomly generated data matrices and other datasets of actual measurements in \cref{section:results}. Finally,  \cref{section:conclusions} concludes the paper. 
%
%
\section{Formulations and Algorithms}
\subsection{Reduced-order modeling, sparse sensing and the previous greedy optimization of sensor placement}
\label{subsection:basics}
First, $p$ observations are linearly constructed from $r_{1}$ parameters as:
\begin{align}
    \bm{y}
    \,&=\,\bm{C}\bm{z}.\label{eq:observation}
\end{align}
Here, $\bm{y}\in \mathbb{R}^{p}$, $\bm{z}\in \mathbb{R}^{r_{1}}$ and $\bm{C}\in \mathbb{R}^{p{\times}r_{1}}$ are an observation vector, a parameter vector and a given measurement matrix, respectively. 
It should also be noted that the absence of noise is assumed in \cref{eq:observation}.
The estimated parameters $\hat{\bm{z}}$ (a quantity with hat refers to the estimated value of the quantity) can be obtained by a pseudo inverse operation.
\begin{align} 
    \hat{\bm{z}}
    \,&=\,\bm{C}^{+}\bm{y}
    \,=\,\left\{\begin{array}{cc}\bm{C}^{\mathrm{T}}\left(\bm{C}\bm{C}^{\mathrm{T}}\right)^{-1}\bm{y}, & p< r_{1}, \\[3pt]
    \left(\bm{C}^{\mathrm{T}}\bm{C}\right)^{-1}\bm{C}^{\mathrm{T}}\bm{y}, &p\geq r_{1}. \end{array} \right.\label{eq:LS_estimation}
\end{align}   
%
Then, uncorrelated Gaussian noise of the same variance and zero mean for every observation point is considered. Joshi and Boyd~\cite{joshi2009sensor} clearly showed an objective function to design the best estimation system for the case $p \geq r_{1}$. They aimed to maximize a logarithm of the determinant of the Fisher information matrix, which realised the least ellipsoid volume of the expected estimation error $\bm{z}-\hat{\bm{z}}$, as in \cref{eq:maximize logdet}
\begin{align}
    \text{maximize}\quad&\log \det \left(\bm{C}^{\mathrm{T}}\bm{C}\right). \label{eq:maximize logdet}
\end{align}
Additionally, there have been many studies adopting those optimization problems to low-order but high-dimensional data. The measurement matrix in \cref{eq:observation} becomes
    $\bm{C}=\bm{H}\bm{C}^{*}$,
where $\bm{C}^{*}$ is a total measurement matrix (or, in other words, a sensor candidate matrix) for whole $n$ observation candidates $(n\gg p,r_{1})$ and $\bm{H}$ as a sensor-location matrix that gives observe locations among $n$ candidates. \Cref{eq:maximize logdet} is now interpreted as a searching problem of the most effective locations of sensors by determining $\bm{H}$ if $\bm{C}^{*}$ is given. Basically, all the combinations of $p$ sensors out of $n$ sensor candidates should be searched by the brute-force algorithm for the real-optimized solution of \cref{eq:maximize logdet}, which takes enormous computational time ($O(n!/(n-p)!/p!)\approx O(n^{p})$).

Instead, greedy algorithms for suboptimized solutions by adding a sensor step by step has been devised for reduced-order modeling~\cite{manohar2018data}. As many studies did, Saito {\it et al.}~\cite{saito2019determinantbased} set $\bm{C}^{*}$ to be the reduced-order spatial-mode matrix $\bm{U}_{1:r_{1}}\in \mathbb{R}^{n{\times}r_{1}}$ in \cref{eq:reduced_order} that proper orthogonal decomposition (POD) generated from given training data, which consists of $m$ snapshots for $n$ variables $\bm{X}\in \mathbb{R}^{n{\times}m}$. With $\bm{U}_{1:r_{1}}$ matrix, a greedy selection was demonstrated as shown in \cref{alg_DGR} to pursue \cref{eq:maximize logdet}.
A data matrix $\bm{X}\in \mathbb{R}^{n{\times}m}\,\left(n>m\right)$ and its reduced-order representation are given by singular value decomposition (SVD):
\begin{align}
    \bm{X}
    \,&=\,\bm{U}\bm{\Sigma}\bm{V}^{\mathrm{T}}\nonumber\\
    \,&=\,\left[\begin{array}{cc}
    \bm{U}_{1:r_{1}}&\bm{U}_{(r_{1}+1):m}
    \end{array}\right]
    \left[\begin{array}{cc}
    \bm{\Sigma}_{1:r_{1}}& \bm{0}\nonumber\\
    \bm{0}&\bm{\Sigma}_{(r_{1}+1):m}
    \end{array}\right]
    \left[\begin{array}{c}
    \bm{V}_{1:r_{1}}^{\mathrm{T}}\nonumber\\\bm{V}_{(r_{1}+1):m}^{\mathrm{T}}
    \end{array}\right]\nonumber\\
    \,&=\,\bm{U}_{1:r_{1}}\bm{\Sigma}_{1:r_{1}}\bm{V}_{1:r_{1}}^{\mathrm{T}}
        +\bm{U}_{(r_{1}+1):m}\bm{\Sigma}_{(r_{1}+1):m}\bm{V}_{(r_{1}+1):m}^{\mathrm{T}}\label{eq:reduced_order}\\
    \,&\equiv\,\bm{X}_{1:r_{1}}+\bm{X}_{(r_{1}+1):m}\nonumber\\
    \,&\approx\,\bm{X}_{1:r_{1}}. \nonumber
\end{align}
$\bm{\Sigma}\in \mathbb{R}^{m{\times}m}$ and $\bm{V}\in \mathbb{R}^{m{\times}m}$ are a diagonal matrix of the singular values and a matrix of the temporal modes, respectively, and the subscript notation $\bm{A}_{i:j}$ for a given matrix $\bm{A}$ denotes a low-rank representation of $\bm{A}$, using the $i$th-to-$j$th singular values or vectors. 
It should be noted that POD can be processed by SVD if spatial and temporal desensitization of the data are uniform.

\begin{algorithm}              
\caption{Overview of DG algorithm~\cite{saito2019determinantbased}}
\begin{algorithmic}
\label{alg_DGR}            

    \STATE$ i_{1}=\argmax_{i\, \in\, \mathcal{S}}
    \bm{u}_i\bm{u}_i^{\mathrm{T}}$
    \STATE$
    \bm{C}_{1}=\bm{u}_{i_{1}}$
\FOR{ $k =2, \dots, r_{1}, \dots, p$ }
\IF{$k \leq r_{1}$}
    \STATE$ i_{k} = \argmax_{i\, \in\, \mathcal{S}\, \backslash\, \mathcal{S}_{k}}\det(\bm{C}_{k}^{\left(i\right)}\bm{C}_{k}^{\left(i\right)\mathrm{T}}) $
      \STATE$ \quad=\argmax_{i\, \in\, \mathcal{S}\, \backslash\, \mathcal{S}_{k}}
        \det\left(\left[ \begin{array}{cc}\bm{C}_{k-1}\\\bm{u}_i\end{array}\right]\left[ \begin{array}{cc} \bm{C}_{k-1}^{\mathrm{T}}& \bm{u}_i^{\mathrm{T}} \end{array} \right] \right)$\\
      \STATE$ \quad= \argmax_{i\, \in\, \mathcal{S}\, \backslash\, \mathcal{S}_{k}}\bm{u}_i\left(\bm{I}-\bm{C}_{k-1}^{\mathrm{T}}\left(\bm{C}_{k-1}\bm{C}_{k-1}^{\mathrm{T}}\right)^{-1}\bm{C}_{k-1}\right)\bm{u}_i^{\mathrm{T}}$
\ELSE
      \STATE$ i_{k}= \argmax_{i\, \in\, \mathcal{S}\, \backslash\, \mathcal{S}_{k}}\det(\bm{C}_{k}^{\left(i\right) \mathrm{T}}\bm{C}_{k}^{\left(i\right)})$
      \STATE$\quad =\argmax_{i\, \in\, \mathcal{S}\, \backslash\, \mathcal{S}_{k}}\det\left(
        \left[ \begin{array}{cc}\bm{C}_{k-1}^{\mathrm{T}}& \bm{u}_i^{\mathrm{T}}\end{array}\right]
        \left[ \begin{array}{cc}\bm{C}_{k-1}\\\bm{u}_i \end{array} \right] \right) $\\
      \STATE$\quad =\argmax_{i\, \in\, \mathcal{S}\, \backslash\, \mathcal{S}_{k}} \left(1 + \bm{u}_i\left(\bm{C}_{k-1}^{\mathrm{T}}\bm{C}_{k-1}\right)^{-1} \bm{u}_i^{\mathrm{T}}\right)$
\ENDIF
    \STATE $ 
    \bm{C}_{k}=\left[\begin{array}{ccc}
        \bm{C}_{k-1}^{\mathrm{T}}&
        \bm{u}_{i_{k}}^{\mathrm{T}}
        \end{array}\right]^{\mathrm{T}}$
\ENDFOR
\end{algorithmic}
\end{algorithm}
Here,
$\mathcal{S}$ and $\mathcal{S}_{k}$ $\left(
    \mathcal{S}=\{1,\hdots,n\}
    ,\,  
    \mathcal{S}_{k}=\{i_{1},\hdots,i_{k}\}
    ,\, 
    \mathcal{S}_{k}\subset\mathcal{S}\right)$
refer a set of indices for locations of the sensor candidates and its subset of the determined sensors, respectively.
Additionally, the notation $\bm{A}_{k}^{\left(i\right)}$ for an arbitrary quantity $\bm{A}$ indicates a quantity that its $k$th component is to be investigated in this step-wise selection, with the $(k-1)$ components given by the $(k-1)$ selected sensors. 
After the $k$th selection step, $\bm{A}_{k}$ is constructed by the $k$ sensors selected.
As an example, $\bm{C}_{k}^{\left(i\right)}$ and $\bm{C}_{k}$ are written as follows:
\begin{eqnarray}
\bm{C}_{k}^{\left(i\right)}
\,&=&\,\left[\begin{array}{ccccc}\bm{u}_{i_{1}}^{\mathrm{T}} &\bm{u}_{i_{2}}^{\mathrm{T}} &\dots &\bm{u}_{i_{k-1}}^{\mathrm{T}} &\bm{u}_{i}^{\mathrm{T}}\end{array}\right]^{\mathrm{T}}\nonumber\\
\bm{C}_{k}
\,&=&\,\bm{C}_{k}^{\left(i\right)}\, |_{i=i_{k}},
\end{eqnarray}
where ${i}_{k}$ and $\bm{u}_{i_{k}} \left(k\in\{1,\hdots,p\}\right)$ indicate an index of the $k$th selected sensor location and the corresponding row vector of the sensor-candidate matrix $\bm{U}_{1:r_{1}}$, respectively. 
In Ref.~\cite{saito2019determinantbased}, maximization of the determinant of the Fisher information matrix of both cases in \cref{eq:LS_estimation} are realized with the matrix determinant lemma which quickens the computation as shown in \cref{alg_DGR}.
%

This approach effectively works, but sometimes does not as experienced in Ref.~\cite{saito2019determinantbased}, and as demonstrated later in a tested problem of \cref{section:results}. That defect arises because the linear least squares estimation that a large number of studies employed, misses information in data such as an expected distribution of parameters $\bm{z}$. Here, $\bm{z}$ in our problem will be amplitudes of principal $r_{1}$ POD modes at arbitrary time. Consequently, the same levels of the amplitudes of low- to high-order modes as each other are assumed although the matrix $\bm{\Sigma}$ in \cref{eq:reduced_order} shows that they differ. Moreover, the aspects are not considered in the previous straightforward implementation that the observation $\bm{y}$ is always contaminated with the truncated POD modes (see matrices in the second term of \cref{eq:reduced_order}), and that they cause spatially correlated noise for the reduced-order estimation \cref{eq:LS_estimation}.
%
\subsection{Bayesian Estimation using Sparse Sensor} \label{subsection:proposed}
In Ref.~\cite{liu2016sensor}, a proper objective function for the selection is derived under an assumption of the correlated noise.
Although they formulated optimization considering such noise, noise itself is only given as a function of distance between sensors.
A data-driven method for the modeling of the correlated noise are presented in this paper and two more conditions are exploited for a more robust estimation; one is expected variance of the POD mode amplitudes, and the other is spatial covariance of the components that are truncated in the order reduction of data matrix \cref{eq:reduced_order}.
The former can be estimated from $\bm{\Sigma}$ as:
\begin{eqnarray}
E(\bm{zz}^{\mathrm{T}})
\,&\equiv&\,\bm{Q}\nonumber\\
\,&\approx&\,\frac{1}{m}\bm{\Sigma}_{1:r_{1}}\bm{V}_{1:r_{1}}^{\mathrm{T}}\bm{V}_{1:r_{1}}\bm{\Sigma}_{1:r_{1}}\nonumber\\
\,&\propto&\,\bm{\Sigma}_{1:r_{1}}^2,
\end{eqnarray}
where $E(\theta)$ is the expectation value of a variable $\theta$.
Then, full-state observation $\bm{x}$ and covariance matrix $\bm{\mathcal{R}}$ of the noise become
\begin{align}
\bm{x}\,&=\,\bm{U}_{1:r_{1}}\bm{z}+\bm{w}\nonumber\\
E(\bm{w}\bm{w}^{\mathrm{T}})\,&\equiv\,\bm{\mathcal{R}}, \label{eq:full cov}
\end{align}
where $\bm{w}$ is a observation noise vector and $\bm{x}$ is one snapshot (one column vector) of $\bm{X}$.
The sparse observation and its noise covariance are:
\begin{eqnarray}
\bm{y}\,&=&\,\bm{H}\bm{U}_{1:r_{1}}\bm{z}+\bm{H}\bm{w},\nonumber\\
E(\bm{H}\bm{w}\bm{w}^{\mathrm{T}}\bm{H}^{\mathrm{T}})\,&=&\,\bm{H}E(\bm{w}\bm{w}^{\mathrm{T}})\bm{H}^{\mathrm{T}}\nonumber\\
\,&\equiv&\,\bm{H}\bm{\mathcal{R}}\bm{H}^{\mathrm{T}}\nonumber\\
\,&\equiv&\,\bm{R}, \label{eq:sparse cov}
\end{eqnarray}
where $\bm{R}\in\mathbb{R}^{p \times p}$  represents a covariance matrix of the noise that $p$ sensors capture.
Here, the full-state noise covariance is assumed to be estimated from the high-order modes: 
\begin{align}
\bm{\mathcal{R}}\,&=\,E(\bm{w}\bm{w}^{\mathrm{T}})\nonumber\\
      \,&=\,E((\bm{x}-\bm{U}_{1:r_{1}}\bm{z})(\bm{x}-\bm{U}_{1:r_{1}}\bm{z})^{\mathrm{T}})\nonumber\\
      \,&\approx\,(\bm{U\Sigma V}^{\mathrm{T}}-\bm{U}_{1:r_{1}}\bm{\Sigma}_{r_{1}}\bm{V}_{r_{1}}^{\mathrm{T}})(\bm{U\Sigma V}^{\mathrm{T}}-\bm{U}_{1:r_{1}}\bm{\Sigma}_{r_{1}}\bm{V}_{r_{1}}^{\mathrm{T}})^{\mathrm{T}}\nonumber\\      \,&=\,(\bm{U}_{(r_{1}+1):m}\bm{\Sigma}_{(r_{1}+1):m}\bm{V}_{(r_{1}+1):m}^{\mathrm{T}})(\bm{U}_{(r_{1}+1):m}\bm{\Sigma}_{(r_{1}+1):m}\bm{V}_{(r_{1}+1):m}^{\mathrm{T}})^{\mathrm{T}}\nonumber\\
      \,&=\,\bm{U}_{(r_{1}+1):m} \bm{\Sigma}^2_{(r_{1}+1):m} \bm{U}_{(r_{1}+1):m}^{\mathrm{T}}, \label{eq:covariance1}      
\end{align}
and
\begin{eqnarray}
\bm{R}\,&=&\,\bm{H}\bm{\mathcal{R}}\bm{H}^{\mathrm{T}}\nonumber\\
\,&\approx&\,\bm{H}(\bm{U}_{(r_{1}+1):m}\bm{\Sigma}^2_{(r_{1}+1):m}\bm{U}_{(r_{1}+1):m}^{\mathrm{T}})\bm{H}^{\mathrm{T}}. \label{eq:covariance2}      
\end{eqnarray}
Then, the Bayesian estimation is derived with those prior information. Here,an a priori probability density function (PDF) of the POD mode amplitudes becomes:
\begin{eqnarray}
P(\bm{z}) \propto \exp (-\bm{z}^{\mathrm{T}}\bm{Q}^{-1}\bm{z}),
\end{eqnarray}
and the conditional PDF of $\bm{y}$ under given $\bm{z}$ is as follows:
\begin{eqnarray}
P(\bm{y}|\bm{z}) \propto \exp (-(\bm{y}-\bm{C}\bm{z})^{\mathrm{T}}\bm{R}^{-1}(\bm{y}-\bm{C}\bm{z})).
\end{eqnarray}
These relations lead to the a posteriori PDF:
\begin{eqnarray}
P(\bm{z}|\bm{y}) &\propto & P(\bm{y}|\bm{z})P(\bm{z})\nonumber\\
&\propto& \exp (-(\bm{y}-\bm{C}\bm{z})^{\mathrm{T}}\bm{R}^{-1}(\bm{y}-\bm{C}\bm{z}))\exp (-\bm{z}^{\mathrm{T}}\bm{Q}^{-1}\bm{z})\nonumber\\
\,&=&\, \exp (-(\bm{y}-\bm{C}\bm{z})^{\mathrm{T}}\bm{R}^{-1}(\bm{y}-\bm{C}\bm{z})-\bm{z}^{\mathrm{T}}\bm{Q}^{-1}\bm{z}).
\end{eqnarray}
Here, the maximum a posteriori estimation on $p(\bm{z}|\bm{y})$ is:
\begin{align} \label{eq:RQ_estimate}
\hat{\bm{z}}=(\bm{C}^{\mathrm{T}}\bm{R}^{-1}\bm{C}+\bm{Q}^{-1})^{-1}\bm{C}^{\mathrm{T}}\bm{R}^{-1}\bm{y}.
\end{align}
Thanks to the normalization term $\bm{Q}$, the inverse operation in \cref{eq:RQ_estimate} is regular for any conditions of $p$ unlike the least squares estimation in \cref{eq:LS_estimation}.
In this estimation, the objective function in the \cref{eq:maximize logdet} is modified:
\begin{align}
\text{maximize} \hspace{10pt} \log\det(\bm{C}^{\mathrm{T}}\bm{R}^{-1}\bm{C}+\bm{Q}^{-1}). \label{eq:detctrc}
\end{align}
Note that sensors should be removed from candidates if they have extremely low signal fluctuations.
This is because those sensors increase the objective value \cref{eq:detctrc} by making the matrix $\bm{R}$ singular. 
In the present study, the locations are beforehand excluded from $\bm{\mathcal{S}}$ for simplicity if their RMSs are $10^2-10^3$ times smaller than the maximum of the dataset, and $\bar{\bm{\mathcal{S}}}$ denotes a subset of $\bm{\mathcal{S}}$ after this exclusion. 

Note that for the case $\bm{R}$ is a diagonal matrix, Joshi and Boyd \cite{joshi2009sensor} have already derived a convex optimization of the sensor selection for the Bayesian estimation. In the present study, $\bm{R}$ includes nondiagonal components which represent the correlation in measurement noise.
%
\begin{algorithm}                        
\caption{Determinant-based greedy algorithm considering noise correlation between sensors}
\begin{algorithmic}
\label{alg_BDG}  
\STATE  $\bm{Q}=\bm{\Sigma}_{1:r_{1}}^2$\\
\STATE $\bm{\mathcal{R}}=\bm{U}_{(r_{1}+1):m} \bm{\Sigma}^2_{(r+1):m} \bm{U}_{(r_{1}+1):m}^{\mathrm{T}}$\\
\FOR{ $k =1, \dots, p$ }
      \STATE $
      i_{k}=\argmax_{i\, \in\, \bar{\mathcal{S}}\, \backslash\, \mathcal{S}_{k}}\det\left(\bm{C}_{k}^{\left( i \right)\mathrm{T}}\left(\bm{R}_{k}^{\left( i \right)}\right)^{-1}\bm{C}_{k}^{\left( i \right)}+\bm{Q}^{-1}\right) 
     \newline
      \quad\left( s.t.\quad \bm{R}_{k}^{\left( i \right)}=\bm{H}_{k}^{\left( i \right)}\bm{\mathcal{R}}\bm{H}_{k}^{\left( i \right)\mathrm{T}}\right)
      $\\
      %
    \STATE $\bm{H}_{k}=\left[\begin{array}{cc}
    \bm{H}_{k-1} \\ \bm{h}_{i_{k}}\end{array}\right], \, 
        \bm{C}_{k}=\left[\begin{array}{ccccc}
            \bm{C}_{k-1}\\[2pt]
            \bm{u}_{i_{k}}
            \end{array}\right]$\\
    \STATE $\bm{R}_{k}=\bm{H}_{k}\bm{\mathcal{R}}\bm{H}_{k}^\mathrm{T}$
\ENDFOR
\end{algorithmic}
\end{algorithm}
The proposed method suppresses correlated measurement noise by considering $\bm{R}$.
This Bayesian determinant-based greedy (BDG) algorithm is presented in \cref{alg_BDG}. Here, $1\times n$ row vector $\bm{h}_{i}$ refers to the $i$th sensor location  that has unity in the $i$th component with zero in the others, which extracts the $i$th row vector from the sensor-candidate matrix $\bm{U}_{1:r_{1}}$.
%
%
%
%
%
%
%
\subsection{Fast algorithm}\label{section:Matrixlemma}
A fast implementation is considered based on \cref{alg_BDG} as Saito {\it et al.} demonstrated in their determinant calculation using rank-one lemma \cite{saito2019determinantbased}.
First, the covariance matrix generated by the $i$th sensor candidate in the $k$th sensor selection, $\bm{R}_{k}^{\left(i\right)}$ is:
\begin{eqnarray}
    \bm{R}_k^{\left(i\right)}\,&=&\,\left(\begin{array}{cc}\bm{R}_{k-1}&\bm{s}_{k}^{\left(i\right)\mathrm{T}}\\\bm{s}_{k}^{\left(i\right)}&{t}_{k}^{\left(i\right)}\end{array}\right),
\end{eqnarray}
where the covariance $\bm{s}_{k}^{\left( i \right)}$ of noise between the sensor candidate $i$ and the $\left(k-1\right)$ selected sensors is
\begin{eqnarray}
\bm{s}_{k}^{\left( i \right)}\,&\propto&\,E(\bm{h}_{i}\bm{ww}^{\mathrm{T}} \bm{H}_{k-1}^{\mathrm{T}}) \nonumber\\
          \,&=&\,\bm{h}_{i}E(\bm{ww}^{\mathrm{T}}) \bm{H}_{k-1}^{\mathrm{T}} \nonumber\\
          \,&=&\,\bm{h}_{i} \bm{\mathcal{R}} \bm{H}_{k-1}^{\mathrm{T}}\nonumber\\
          \,&\approx&\,\bm{h}_{i}(\bm{U}_{(r_{1}+1):m}\bm{\Sigma}^2_{(r+1):m}\bm{U}_{(r_{1}+1):m}^{\mathrm{T}})\bm{H}_{k-1}^{\mathrm{T}} \label{eq:p of R},
\end{eqnarray}
and similarly, the variance of noise at the location $i$ is
\begin{eqnarray}
{t}_{k}^{\left( i \right)}
       \,&\approx&\,\bm{h}_{i}(\bm{U}_{(r_{1}+1):m}\bm{\Sigma}^2_{(r+1):m}\bm{U}_{(r_{1}+1):m}^{\mathrm{T}})\bm{h}_{i}^{\mathrm{T}}. \label{eq:q of R}
\end{eqnarray}
Here, $\bm{H}_{k-1}$ is the optimized first-to-($k-1$)th-sensor selection matrix.
In addition, $\bm{R}_{k-1}$ is determined in the previous $\left(k-1\right)$th step. Accordingly, $\left(\bm{R}_{k}^{\left( i \right)}\right)^{-1}$ is obtained as follows:
\begin{eqnarray}
    \left(\bm{R}_{k}^{\left( i \right)}\right)^{-1}
    \equiv \left(\begin{array}{cc}\bm{\alpha}_{k}^{\left( i \right)}&\bm{\beta}_{k}^{\left( i \right)\mathrm{T}}\\[1pt]\bm{\beta}_{k}^{\left( i \right)}&{\delta}_k^{\left( i \right)}\end{array}\right),
\end{eqnarray}
where
\begin{eqnarray}
    \bm{\alpha}_{k}^{\left( i \right)}\,&=
    &\,\bm{R}_{k-1}^{-1}+\frac{1}{{t}_{k}^{\left( i \right)}-\bm{s}_{k}^{\left( i \right)}\bm{R}_{k-1}^{-1}\bm{s}_{k}^{\left( i \right)\mathrm{T}}}\bm{R}_{k-1}^{-1}\bm{s}_{k}^{\left( i \right)\mathrm{T}}\bm{s}_{k}^{\left( i \right)}\bm{R}_{k-1}^{-1},\nonumber\\
    \bm{\beta}_{k}^{\left( i \right)}\,&=
    &\,-\frac{\bm{s}_{k}^{\left( i \right)}\bm{R}_{k-1}^{-1}}{{t}_{k}^{\left( i \right)}-\bm{s}_{k}^{\left( i \right)}\bm{R}_{k-1}^{-1}\bm{s}_{k}^{\left( i \right)\mathrm{T}}},\nonumber\\
    {\delta}_{k}^{\left( i \right)}\,&=
    &\,\frac{1}{{t}_{k}^{\left( i \right)}-\bm{s}_{k}^{\left( i \right)}\bm{R}_{k-1}^{-1}\bm{s}_{k}^{\left( i \right)\mathrm{T}}}.\nonumber
\end{eqnarray}
The objective function is now considered based on the expressions above.
\begin{align}
     &\det\left(\bm{W}_{k}^{\left(i\right)}\right)\nonumber\\
     \,&\equiv\,\det{\left(\bm{C}_{k}^{\left(i\right)\mathrm{T}}\left(\bm{R}_{k}^{\left( i \right)}\right)^{-1}\bm{C}_{k}^{\left(i\right)}+\bm{Q}^{-1}\right)}\nonumber\\
     \,&=\,\det\left(\left[\begin{array}{cc}\bm{C}_{k-1}^{\mathrm{T}}&\bm{u}_{i}^{\mathrm{T}}\end{array}\right]\left[\begin{array}{cc}\bm{\alpha}_{k}^{\left(i\right)}&\bm{\beta}_{k}^{\left(i\right)\mathrm{T}}
     \\ \bm{\beta}_{k}^{\left(i\right)}&{\delta}_{k}^{\left(i\right)}\end{array}\right]\left[\begin{array}{c}\bm{C}_{k-1}\\\bm{u}_{i}\end{array}\right]+\bm{Q}^{-1}\right)\nonumber\\
     \,&=\,\det\left(\bm{C}_{k-1}^{\mathrm{T}}\bm{\alpha}_{k}^{\left(i\right)}\bm{C}_{k-1}+\bm{u}_{i}^{\mathrm{T}}\bm{\beta}_{k}^{\left(i\right)}\bm{C}_{k-1}+\bm{C}_{k-1}^{\mathrm{T}}\bm{\beta}_{k}^{\left(i\right)\mathrm{T}}\bm{u}_{i}+{\delta}_{k}^{\left(i\right)}\bm{u}_{i}^{\mathrm{T}}\bm{u}_{i}+\bm{Q}^{-1}\right)\nonumber\\
     \,&=\,\det\left(\bm{W}_{k-1}+\frac{\left(\bm{C}_{k-1}^{\mathrm{T}}\bm{R}_{k-1}^{-1}\bm{s}_{k}^{\left(i\right)\mathrm{T}}-\bm{u}_{i}^{\mathrm{T}}\right)\left(\bm{s}_{k}^{\left(i\right)}\bm{R}_{k-1}^{-1}\bm{C}_{k-1}-\bm{u}_{i}\right)}{{t}_{k}^{\left(i\right)}-\bm{s}_{k}^{\left(i\right)}\bm{R}_{k-1}^{-1}\bm{s}_{k}^{\left(i\right)\mathrm{T}}}\right)\nonumber\\
     \,&=\,\left(1+
     \frac
     {\left(\bm{s}_{k}^{\left(i\right)}\bm{R}_{k-1 }^{-1}\bm{C}_{k-1}-\bm{u}_{i}\right)
        \bm{W}_{k-1}^{-1}
        \left(\bm{C}_{k-1}^{\mathrm{T}}\bm{R}_{k-1}^{-1}\bm{s}_{k}^{\left(i\right)\mathrm{T}}-\bm{u}_{i}^{\mathrm{T}}\right)}
     {{t}_{k}^{\left(i\right)}-\bm{s}_{k}^{\left(i\right)}\bm{R}_{k-1}^{-1}\bm{s}_{k}^{\left(i\right)\mathrm{T}}}\right)\\
    & \hspace{160pt}\times\det\left(\bm{W}_{k-1}\right).\nonumber
\end{align}
Because $\bm{W}_{k-1}$ and its determinant have already been obtained in the previous step, the $k$th sensor can be selected by the following scalar evaluation:
\begin{align}\label{eq:RQ_sensor} 
     \argmax_{i\, \in\, \bar{\mathcal{S}}\backslash\mathcal{S}_{k}}
     \frac{
        \left(\bm{s}_{k}^{\left(i\right)}\bm{R}_{k-1}^{-1}\bm{C}_{k-1}-\bm{u}_{i}\right)\bm{W}_{k-1}^{-1}\left(\bm{C}_{k-1}^{\mathrm{T}}\bm{R}_{k-1}^{-1}\bm{s}_{k}^{\left(i\right)\mathrm{T}}-\bm{u}_{i}^{\mathrm{T}}\right)
        }{
        {t}_{k}^{\left(i\right)}-\bm{s}_{k}^{\left(i\right)}\bm{R}_{k-1}^{-1}\bm{s}_{k}^{\left(i\right)\mathrm{T}}
        }.
\end{align}
Once the sensor is selected, $\bm{R}_{k}^{-1}$ and $\bm{W}_{k}$ are updated. This algorithm is described in \cref{alg_fastBDG} including another improvement for computational efficiency which is introduced in the next subsection.

\subsection{Memory Efficient Implementation} \label{section:r3truncation}
As was already introduced in \cref{eq:covariance1}, the covariance of noise for every pair of observation points can be estimated by multiplying $\bm{U}_{(r_{1}+1):m}$ and $\bm{\Sigma}_{(r_{1}+1):m}$, and then stored to be $\bm{\mathcal{R}}$. In \cref{alg_BDG}, $\bm{s}_{k}^{\left(i\right)}$ and ${t}_{k}^{\left(i\right)}$ are constructed by taking the corresponding parts of $\bm{\mathcal{R}}$. 
Although this could be a straightforward way to calculate $\bm{s}_{k}^{\left(i\right)}$ and ${t}_{k}^{\left(i\right)}$, it often runs out memory storage to store $\bm{\mathcal{R}}$ since $n^{2}$, the size of $\bm{\mathcal{R}}$, can often reach billions or more in an actual application.

Therefore, the following implementation is proposed by reducing the order of the noise covariance from $\left(m-r_{1}\right)$ to $r_{2}$: only diagonal components of $\mathcal{\bm{R}}$ is constructed and stored as a $n$-components vector $\bm{d}$, and $\bm{s}_{k}^{\left(i\right)}$ is approximated for every $i$ loop by using the dominant first $r_{2}$ columns of $\bm{U}_{(r_{1}+1):m}\bm{\Sigma}_{(r_{1}+1):m}$ in \cref{eq:covariance1}. This modification is reasonable because nondiagonal components are inherently small compared to diagonal ones of $\mathcal{\bm{R}}$, thereby modes have less effect on the $\bm{R}_{k}$ as the mode number increases.
Here, an appropriate $r_{2}$ should be determined with consideration on characteristics of the data matrix.
The effects of truncating $r_{2}$ of $\bm{s}_{k}^{\left(i\right)}$ are shown in \cref{section:sst}.

Then, $\bm{s}_{k}^{\left(i\right)}$ is approximated:
\begin{eqnarray}\label{eq:p of R trunc}
    \tilde{\bm{s}_{k}}^{\left(i\right)}
    \,\approx\,(\bm{h}_{i}\bm{U}_{(r_{1}+1):(r_{1}+r_{2})})\bm{\Sigma}^2_{(r_{1}+1):(r_{1}+r_{2})}(\bm{U}_{(r_{1}+1):(r_{1}+r_{2})}^{\mathrm{T}}\bm{H}_{k-1}^{\mathrm{T}}),
\end{eqnarray}
where, $\bm{U}_{(r_{1}+1):(r_{1}+r_{2})}$ and $\bm{\Sigma}_{(r_{1}+1):(r_{1}+r_{2})}$ are the leading $r_{2}$ columns of the remainder spatial modes matrix $\bm{U}_{(r_{1}+1):m}$ and first $r_{2}\times r_{2}$ components of the remainder singular value matrix $\bm{\Sigma}_{(r_{1}+1):m}$.
These processes are called the $r_{2}$ truncation. This approximation simply reduces amount of stored memory and computational complexity. 
These modifications in \cref{section:Matrixlemma} and \cref{section:r3truncation} are integrated into the \cref{alg_fastBDG}.
A comparison on computing complexity of the methods introduced thus far are listed in \cref{table:sensor_selection_comlexity}.
\begin{algorithm}                        
\caption{Detailed accelerated determinant-based greedy algorithm considering correlation between sensors}
\begin{algorithmic}
\label{alg_fastBDG}  
\STATE Set amplitudes variance matrix \\
\STATE \qquad $\bm{Q}=\bm{\Sigma}_{1:r_{1}}^2$\\
\STATE Set noise variance vector \\
\STATE \qquad
$\bm{d}\quad\left( s.t.\, \bm{d}(j)=\bm{h}_{j}\bm{U}_{(r_{1}+1):m}\bm{\Sigma}^2_{(r+1):m}\bm{U}_{(r_{1}+1):m}^{\mathrm{T}}\bm{h}_{j}^{\mathrm{T}}\right)$\\
\STATE $i_{1}=\argmax_{i\, \in\, \bar{\mathcal{S}}}\det\left(\bm{u}_{i}^{\mathrm{T}}{t}_{1}^{\left(i\right)-1}{\bm{u}}_{i}+\bm{Q}^{-1}\right) 
      \newline \hspace*{8pt}
      = \argmax_{i\, \in\, \bar{\mathcal{S}}} {\bm{u}_{i}\bm{Q}\bm{u}_{i}^{\mathrm{T}}}/{{t}_{1}^{\left(i\right)}}$\\
       \STATE$\quad\left( s.t. \quad
      \bm{t}_{1}^{\left(i\right)} =\bm{h}_{i}\bm{d}^{\mathrm{T}}
      \right)$\\
\STATE Set sensor-location and observation matrix \\[2pt]
\STATE \qquad$\bm{H}_{1}=\bm{h}_{i_{1}},
    \,
    \bm{C}_{1}=\bm{u}_{i_{1}}$\\ 
\STATE Set sensor-covariance matrix \\[2pt]
\STATE \qquad$\bm{R}_{1}=      \bm{h}_{i_{1}}\bm{d}^{\mathrm{T}}$\\
\FOR{ $k =2, \dots, r, \dots, p$ }
    \STATE Calculate and store $(\bm{C}_{k-1}^{\mathrm{T}}\bm{R}_{k-1}^{-1}\bm{C}_{k-1}+\bm{Q}^{-1})^{-1}$, $\bm{R}_{k-1}^{-1}\bm{C}_{k-1}$
      \STATE $i_{k}=\argmax_{i\, \in\, \bar{\mathcal{S}}\, \backslash\, \mathcal{S}_{k}}\det\left(\bm{C}_{k}^{\mathrm{T}}  \left(\bm{R}_{k}^{\left(i\right)}\right)^{-1}\bm{C}_{k} +\bm{Q}^{-1}\right) 
      \newline \hspace*{8pt}
      = \argmax_{i\, \in\, \bar{\mathcal{S}}\, \backslash\, \mathcal{S}_{k}} \left(\bm{s}_{k}^{\left(i\right)}\bm{R}_{k-1}^{-1}\bm{C}_{k-1}-\bm{u}_{i}\right)\left(\bm{C}_{k-1}^{\mathrm{T}}\bm{R}_{k-1}^{-1}\bm{C}_{k-1}+\bm{Q}^{-1}\right)^{-1}
      \newline \hspace*{90pt}
      \times\left(\bm{C}_{k-1}^{\mathrm{T}}\bm{R}_{k-1}^{-1}\bm{s}_{k}^{\left(i\right)\mathrm{T}}-\bm{u}_{i}^{\mathrm{T}}\right)/\left({t}_{k}^{\left(i\right)}-\bm{s}_{k}^{\left(i\right)}\bm{R}_{k-1}^{-1}\bm{s}_{k}^{\left(i\right)\mathrm{T}}\right)$
        \newline
      $\quad\left(
      \begin{array}{rl}
      s.t. 
      \,\,
      \bm{s}_{k}^{\left(i\right)} =&\left(\bm{h}_{i}\bm{U}_{(r_{1}+1):(r_{1}+r_{2})}\right)\bm{\Sigma}^2_{(r_{1}+1):(r_{1}+r_{2})}\left(\bm{U}_{(r_{1}+1):(r_{1}+r_{2})}^{\mathrm{T}}\bm{H}_{k-1}^{\mathrm{T}}\right)
      \\
      \bm{t}_{k}^{\left(i\right)} =&\bm{h}_{i}\bm{d}^{\mathrm{T}}
      
      \end{array}
      \right)$\\

\STATE Set sensor-location and observation matrix \\[2pt]
\STATE \qquad$
    \bm{H}_{k}=\left[\begin{array}{ccccc}
    \bm{H}_{k-1}\\[2pt]
    \bm{h}_{i_{k}}
    \end{array}\right]
    ,\,
    \bm{C}_{k}=\left[\begin{array}{ccccc}
    \bm{C}_{k-1}\\[2pt]
    \bm{u}_{i_{k}}
    \end{array}\right]$
\STATE Set noise-covariance matrix \\[2pt]
\STATE \qquad$\begin{array}{rl}
    \bm{s}_{k}=&\left(\bm{h}_{i_{k}}\bm{U}_{(r_{1}+1):(r_{1}+r_{2})}\right)\bm{\Sigma}^2_{(r_{1}+1):(r_{1}+r_{2})}\left(\bm{U}_{(r_{1}+1):(r_{1}+r_{2})}^{\mathrm{T}}\bm{H}_{k-1}^{\mathrm{T}}\right)\\
    {t}_{k}=&\bm{h}_{i_{k}}\bm{d}^{\mathrm{T}}\\
    \bm{R}_{k}=& \left(\begin{array}{cc}
        \bm{R}_{k-1} 
        & 
        \bm{s}_k^{\mathrm{T}}
        \\
        \bm{s}_k
        &
        {t}_k
        \end{array}\right)
    \end{array}$ \\
\ENDFOR
\end{algorithmic}
\end{algorithm}
\begin{table}[htbp]
\small
\begin{center}
\caption{Computational complexity for four different sensor selection methods} 
\label{table:sensor_selection_comlexity}
\setlength{\tabcolsep}{3pt}
\begin{tabular}{|p{60pt}|p{90pt}p{70pt}|}
\hline
\textbf{Name}
&
\textbf{Complexity}
&
\\
\hline\hline
Brute-force
&
${\frac{n!}{\left(n-p\right)!p!}} \sim {O\left(n^p\right)}$
&
\\
\hline
DG&
$\begin{array}{l}
    {p \le r_{1}} : {O\left(np^2\right)}
    \\
    {p>r_{1}} : {O\left(npr_{1}^2\right)}
\end{array}$
&[\cref{alg_DGR}]
\\

\hline

BDG
&
${O\left(np^3r_{1}\right)}$
&[\cref{alg_BDG}]
\\

\hline

Fast-BDG
&
${O\left(np^3\right)}$
&
[\cref{alg_fastBDG}]
\\
\hline
\end{tabular}
\end{center}
\end{table}

\section{Applications} \label{section:results}
\subsection{Sensors for randomized matrix}
\label{section:random}
The numerical experiments were conducted and the proposed method were verified. The random data matrices, $\bm{X}_{\text{rand}}=\bm{U}\bm{\Sigma}\bm{V}^{\mathrm{T}} (\bm{X}_{\text{rand}}\in\mathbb{R}^{1000 \times 500})$, were set, where $\bm{U}$ and $\bm{V}$ consist of 500 orthogonal vectors that were generated by QR decomposition of normally distributed random matrices, and components of diagonal matrix $\bm{\Sigma}$ are $\text{diag}(\bm{\Sigma})=[1,1/\sqrt{2},1/\sqrt{3},...,1/\sqrt{500}]$, respectively. These slowly decaying diagonal components of $S$ are simulating the data of an actual flow field. 
The following results of measuring computational time to obtain sensors demonstrated in \cref{section:random,section:sst} are conducted under the environments as listed in \cref{table:Computational enviroment}. 
\begin{center}
\begin{threeparttable}[htbp]
\centering
\small
\caption{Computing environments to measure computational time} 
\label{table:Computational enviroment}
\begin{tabular}{|p{100pt}|p{95pt}|p{95pt}|}
\hline
\textbf{Specification} & \textbf{Randomized matrix} & \textbf{NOAA-SST}\\
\hline\hline
 Processor information & Intel(R) Core(TM) & Intel(R) Core(TM)\\
  & $i7-2600$S@$2.80$GHz &$i7-6800$K@$3.40$GHz\\
\hline
Random access memory & $4$ GB &  $128$ GB \\
\hline
System type & $64$ bit operating system  & $64$ bit operating system \\
& x64 base processor  &x64 base processor \\
\hline
Program code & \multicolumn{2}{c|}{Matlab R2013a}\\
\hline
Operating system & Linux Mint Tessa & Windows 10 Pro \\
&  Version: 19.1 &Version:1890 \\
\hline
\end{tabular}  
\end{threeparttable}
\end{center}
First, the reconstruction error which is used throughout the paper is defined as follows:
\begin{align}
    \epsilon
    \,&=\,\sum^{m}_{j=1}\frac{\| \bm{x}(j)-\tilde{\bm{x}}(j)\|_{2}}{\| \bm{x}(j)\|_{2}}, \label{eq:error}
\end{align}
where the numerator and denominator are $L_2$ norm of residual and that of the original state, respectively.

\begin{figure}[htbp]
    \centering
    \subfloat[Computational time] {\label{fig:rand_time}\includegraphics[height=2in]{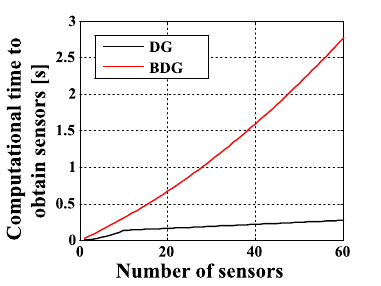}}\\
    \subfloat[Error]{\label{fig:rand_error}\includegraphics[height=2in]{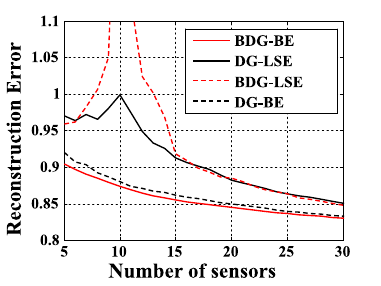}}
    \caption{Sensor selection results on randomized matrix}
\end{figure}

\Cref{fig:rand_time} illustrates the computational time for the sensor selection. This figure shows the increases in computational cost as indicated in  \cref{table:sensor_selection_comlexity}.
\begin{table}[htbp]
\centering
\caption{Collection of introduced methods} 
\setlength{\tabcolsep}{3pt}
\subfloat[State estimation methods]{\label{table:state_estimation_methods}\begin{tabular}{|p{25pt}|p{160pt}p{65pt}|}
\hline
LSE
&
\,$\hat{\bm{z}}=\bm{C}^{+}\bm{y}$&\\
&
    \quad$=\left\{
        \begin{array}{ll}\bm{C}^{\mathrm{T}}\left(\bm{C}\bm{C}^{\mathrm{T}}\right)^{-1}\bm{y} & p\le r \\[3pt]
    \left(\bm{C}^{\mathrm{T}}\bm{C}\right)^{-1}\bm{C}^{\mathrm{T}}\bm{y} &p>r\end{array} \right.$
    &
    [\cref{eq:LS_estimation}]\\

\hline
BE
&
\,$\hat{\bm{z}}=
    \left(\bm{C}^{\mathrm{T}}\bm{R}^{-1}\bm{C}+\bm{Q}^{-1}\right)^{-1}\bm{C}^{\mathrm{T}}\bm{R}^{-1}\bm{y}$&
    [\cref{eq:RQ_estimate}]
\\
\hline
\end{tabular}}\\
\subfloat[Sensor selection methods]{\label{table:sensor_selection_methods}
\begin{tabular}{|p{25pt}|p{160pt}p{65pt}|}
\hline
DG
&
$\begin{array}{ll}
    &\\[-20pt]
    {\argmax\det\left(\bm{CC}^{\mathrm{T}}\right)}
    & p \le r
    \\
    {\argmax\det\left(\bm{C}^{\mathrm{T}}\bm{C}\right)}
    & p>r
    \end{array}$
& 
[\cref{alg_DGR}]
\\
\hline
BDG
&
$\begin{array}{ll}
    &\\[-20pt]
    {\argmax\det\left(\bm{C}^{\mathrm{T}}\bm{R}^{-1}\bm{C}+\bm{Q}^{-1}\right)}\end{array}$
    &
    [\cref{alg_BDG}]
\\
\hline
\end{tabular}}
\end{table}

In addition, the comparison on estimation error defined in \cref{eq:error} was conducted. 
The effectiveness of Bayesian estimation in \cref{eq:RQ_estimate} and the BDG sensor selection in \cref{eq:RQ_sensor} were independently investigated in \cref{fig:rand_error}. 
The estimation methods and the sensor selection methods are presented in \cref{table:state_estimation_methods,table:sensor_selection_methods}, and note that DG-LSE (DG and the least squares estimation) and BDG-BE (BDG and the Bayesian estimation) are optimum combinations under the framework of LSE and BE, respectively.
Two more combinations of the methods were tested for verification, namely DG-BE and BDG-LSE. 

%
%
\Cref{fig:rand_error} illustrates that the Bayesian estimation plays an important role, since the DG-BE result is better than that of DG-LSE although correlated noise is not considered in the selection.
In contrast, the error of LSE increases around $p=10\left(=r_{1}\right)$ regardless of the choice of DG or BDG. 
Errors around $p=r_{1}$ increase because the observed signal by the sensors is strictly converted into the latent amplitudes of $r_{1}$ POD modes despite they contain intense correlated noise of higher modes.
The components in the smallest-sensitivity direction should also be estimated which includes the larger error due to the noise, whereas the components in that direction are assumed to be zero in the pseudo-inverse operation and such an error is suppressed when $p < r_{1}$. 
The error of BDG-LSE is the largest in the range of $6<p<17$ because the sensors of BDG are chosen with assuming the regularization of the Bayesian estimation, and therefore the smallest-sensitivity direction components of least squares estimation are not accurately predicted in BDG-LSE than in DG-LSE. 
However, BDG-LSE works sligtly better than, or at least equal to, DG-LSE in the condition of $p \ge 17$. This is because BDG choose the sensor which is less contaminated by correlated noise and such sensors work better even in the LSE. 
Finally, the error of BDG-BE is always the smallest in these four methods. This illustrates that BDG choose the sensor positions suitable to the Bayesian estimation \cref{eq:RQ_estimate}.

\subsection{Sensors for flow around airfoil}\label{section:piv}
The particle image velocimetry (PIV) was previously conducted and time-resolved data of velocity fields around an airfoil were obtained\cite{nankai2019linear}. The effectiveness of the present method for the PIV data is demonstrated hereafter. The test conditions are listed in \cref{table:PIV_test_conditions} and the fluctuating components of the freestream direction velocity is only employed, unlike the previous study in which the two-dimensional velocity is simultaneously treated \cite{saito2020data}. 

\begin{table}[htbp]
\centering
\caption{PIV test conditions\cite{nankai2019linear}} 
\label{table:PIV_test_conditions}
\setlength{\tabcolsep}{3pt}
\begin{tabular}{|p{150pt}|p{90pt}|}
\hline
Laser& Double pulse lasers\\\hline
Time between pulse& $100$ $\rm{\mu}$s\\\hline
Sampling rate& $5000$ Hz\\\hline
Particle image resolution& $1024 \times 1024$ pixels\\\hline
Total number of image pairs& $9700$\\
\hline
\end{tabular}
\end{table}

\begin{figure}[hbtp]
    \centering
    \subfloat[Position of 20 sensors on RMS map; (1) DG, (2) BDG]{\label{fig:PIV_RMS}\includegraphics[height=2in]{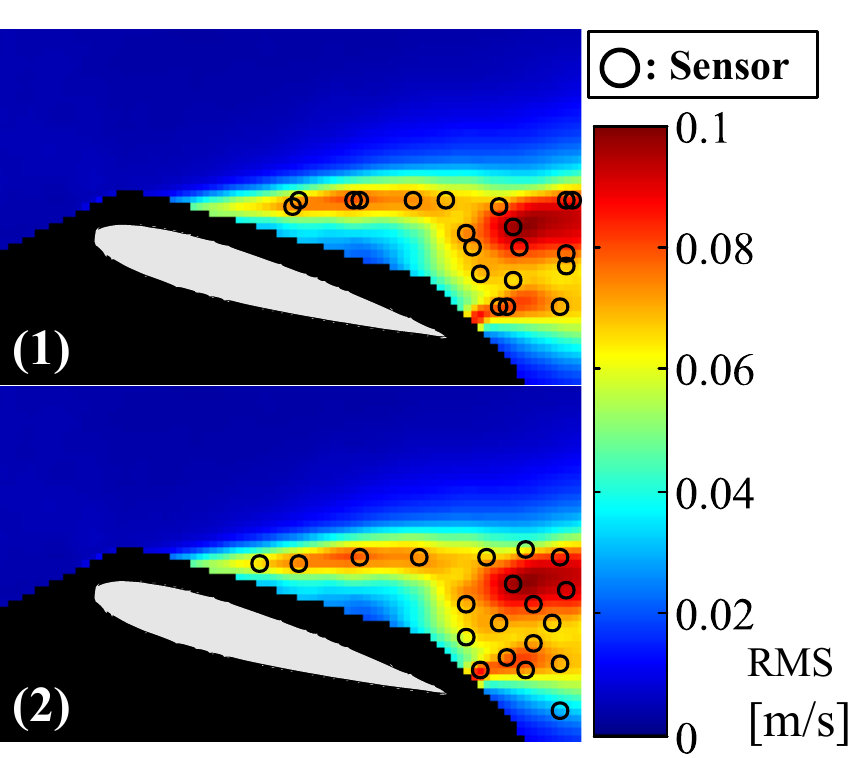}}\\
    \subfloat[Error]{\label{fig:PIV_Error}\includegraphics[height=2in]{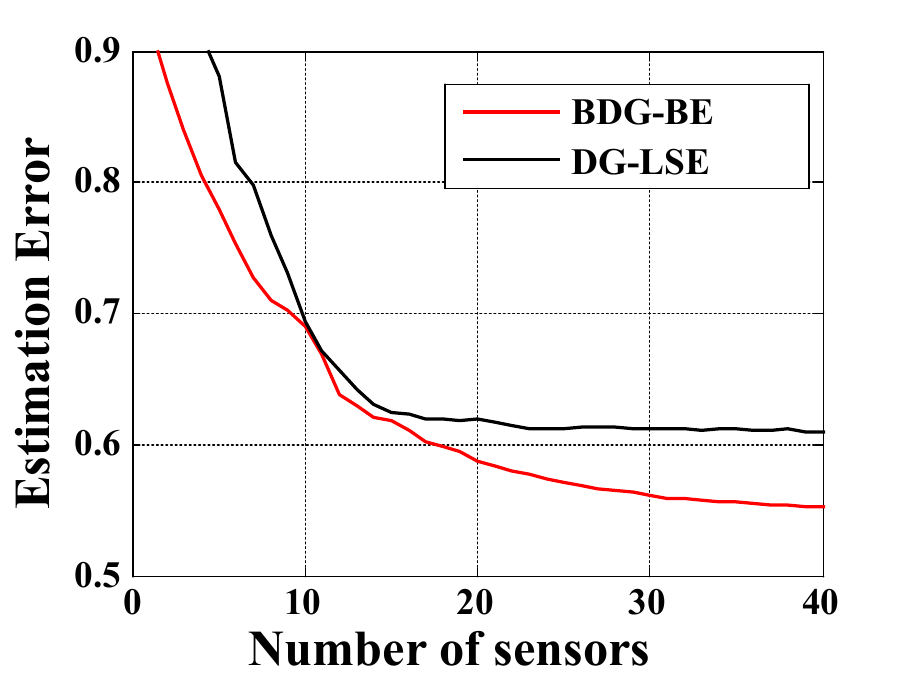}}
    \caption{Sensing of flow around airfoil}
\end{figure}


\Cref{fig:PIV_RMS} clearly shows that the present method selects different sensor positions from those selected by the previous algorithm.
The region of high RMS in the temporal series data is red-colored in the figure.
The selected sensor positions of the previous method tend to be concentrated in small parts of highly fluctuated recirculation regions, whereas those of the presented method are located evenly in those regions.
A comparison of the reconstruction error using the first 10 modes illustrates that the proposed method is effective and its error approaches the lower limit of error estimated by full observation.

\subsection{Sensors for sea surface temperature distribution}
\label{section:sst}
The sea surface temperature (SST) data are distributed at the NOAA website \cite{noaa}. The data formatted by Manohar {\it et al.} \cite{manohar2018data} are adopted in the present study.
The details of the data are listed in \cref{table:SST_data_conditions}.

\begin{table}[htbp]
\centering
\caption{SST data conditions} 
\label{table:SST_data_conditions}
\setlength{\tabcolsep}{3pt}
\begin{tabular}{|p{90pt}|p{220pt}|}
\hline
Brief Description& NOAA Optimum Interpolation (OI) SST V2 \cite{noaa}
\\\hline
Temporal Coverage& Weekly means from 1989/12/31 to 1999/12/12 \\
&
\qquad\qquad\qquad\qquad\qquad\qquad\qquad ($520$ snapshots)\\\hline
Spatial Coverage& 
1.0 degree latitude x 1.0 degree longitude\\
&
\qquad\qquad global grid ($n=44219$ observed points)\\
\hline
\end{tabular}
\end{table}
%

\begin{figure}[htbp]
    \centering
    \subfloat[Position of 20 sensors on RMS map; (1) DG, (2) BDG] {\label{fig:SST_RMS}\includegraphics[clip,height=2in]{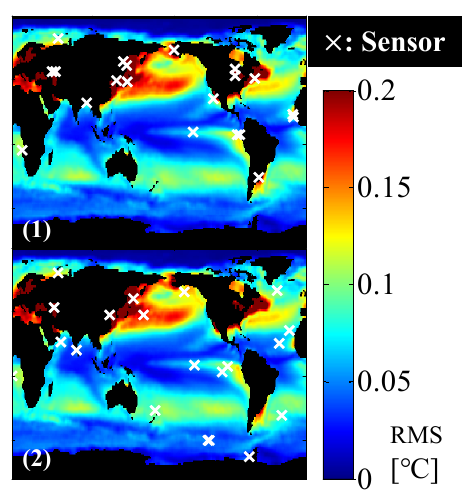}}\\
    \subfloat[Error] {\label{fig:SST_Error}\includegraphics[clip,height=2in]{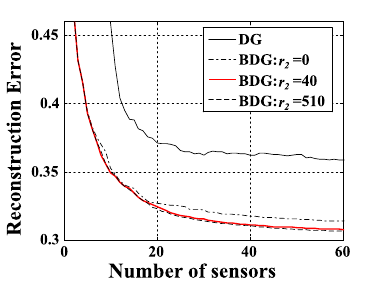}}\\
    \subfloat[Computational time] {\label{fig:SST_r3trunc}\includegraphics[clip,height=2in]{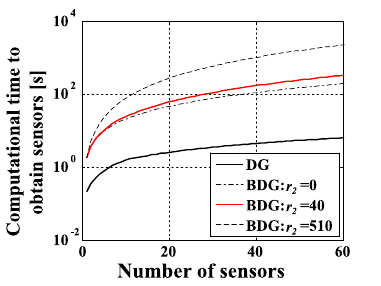}}
    \caption{Sensing of global sea surface temperature}
\end{figure}

As was conducted in \cref{section:piv}, ten modes are employed for the reduced-order modeling after trimming the low RMS points off  from the sensor candidates $\mathcal{S}$ as described in \cref{subsection:proposed}.
\cref{fig:SST_RMS} shows the difference in sensor positions as indicated by cross marks, and backgrounds are time-series RMS of SST data.
Those figures show that the sensors are so scattered by the proposed method that the effects of the correlated noise due to truncated POD modes could be minimized. 

The effect of $r_{2}$ truncation for an efficient memory implementation is also demonstrated in this section as introduced in \cref{section:r3truncation}.
Because the SST dataset consists of a enormous number of observation locations $O\left(10^{4}\right)$, and it requires considerable time to calculate \cref{eq:full cov}. 
The variation of the reconstruction error and the computational time are shown in \cref{fig:SST_Error,fig:SST_r3trunc}.
Those results of the BDG sensors are calculated for three cases, $r_{2}=0$, $40$, $510$. Especially, $r_{2}=0$ and $r_{2}=510$ correspond to the case ignoring all the nondiagonal terms of $\bm{R}_{k}$, and the case  \cref{eq:full cov}, respectively.
As expected, these plots show that more accurate estimation is realized when higher modes are used for noise covariance, but it indeed needs longer calculation time.
A reasonable $r_{2}$ which leads to both less calculation time and less reconstruction error appears to be around $r_{2}=40$, only $10\%$ of the original modes.
The general criteria of reasonable $r_{2}$ parameter should be addressed in future research. 

\section{Conclusions}\label{section:conclusions}
In the context of Bayesian maximum a posteriori (MAP) state estimation, a greedy selection method for noise-robust sensors is proposed for the reconstruction of the high-dimensional data.
\reviewerA{The algorithm leverages high-order modes of the singular value decomposition of data, which were ignored in the previous reduced-order sensing framework, for the construction of the noise covariance matrix used as weights in the MAP state estimation.}
\reviewerA{Prior distribution of coefficients are also generated from information of the modes.}
The proposed method determines sensors by maximizing the determinant of a matrix in the MAP estimation operator, which simultaneously minimizes a metric of expected estimation error.
\reviewerA{Additionally, the rank-one lemma and a low-rank approximation of the covariance matrix are considered for more efficient implementation.}

\reviewerA{Our approach was applied to some examples and the reconstruction ability of determined sensors were verified. Firstly, we separately confirmed noise-robustness of sensors of the proposed method and benefits of the Bayesian estimation. 
Reproduction results on data matrices showed significant improvement in reconstruction accuracy with reasonable increase in the computational time, compared with that of the previously presented method.}

\reviewerA{
Future works of interest are including;
\begin{itemize}
    \item To extend the method to nonlinear measurement.
    This method is based on a linear modeling of POD modes, hence it should be extended to more generalized measurement configuration like~\cite{chaturantabut2010nonlinear}.
    \item To adopt various noise models such as an exponential model from the same reason above.
    \item To involve dynamics of phenomena into modeling, estimation and sensor selection.
    This should be considered for feedback control. Some metrics associated with observability and controllablity were proposed in Ref.~\cite{summers2015submodularity}, and the extension of the present work to these metrics is in our future focus.
\end{itemize}
}

\section*{Acknowledgements}
This work was partially supported by JST CREST Grant Number JPMJCR1763, Japan.
The fourth author T.N. is grateful for support of the grant JPMJPR1678 of JST Presto, Japan.


\bibliography{xaerolab_1016}

\begin{thebibliography}{10}
\expandafter\ifx\csname url\endcsname\relax
  \def\url#1{\texttt{#1}}\fi
\expandafter\ifx\csname urlprefix\endcsname\relax\def\urlprefix{URL }\fi
\expandafter\ifx\csname href\endcsname\relax
  \def\href#1#2{#2} \def\path#1{#1}\fi

\bibitem{kincaid2002d}
R.~K. Kincaid, S.~L. Padula, D-optimal designs for sensor and actuator
  locations, Computers \& Operations Research 29~(6) (2002) 701--713.
\newblock \href {https://doi.org/10.1016/S0305-0548(01)00048-X}
  {\path{doi:10.1016/S0305-0548(01)00048-X}}.

\bibitem{corke2007sdbd}
T.~C. Corke, M.~L. Post, D.~M. Orlov, Sdbd plasma enhanced aerodynamics:
  concepts, optimization and applications, Progress in Aerospace Sciences
  43~(7-8) (2007) 193--217.
\newblock \href {https://doi.org/10.1016/j.paerosci.2007.06.001}
  {\path{doi:10.1016/j.paerosci.2007.06.001}}.

\bibitem{bates1996experimental}
R.~Bates, R.~Buck, E.~Riccomagno, H.~Wynn, Experimental design and observation
  for large systems, Journal of the Royal Statistical Society: Series B
  (Methodological) 58~(1) (1996) 77--94.
\newblock \href {https://doi.org/10.1111/j.2517-6161.1996.tb02068.x}
  {\path{doi:10.1111/j.2517-6161.1996.tb02068.x}}.

\bibitem{brunton2014compressive}
S.~L. Brunton, J.~H. Tu, I.~Bright, J.~N. Kutz, Compressive sensing and
  low-rank libraries for classification of bifurcation regimes in nonlinear
  dynamical systems, SIAM Journal on Applied Dynamical Systems 13~(4) (2014)
  1716--1732.

\bibitem{udwadia1994methodology}
F.~E. Udwadia, Methodology for optimum sensor locations for parameter
  identification in dynamic systems, Journal of engineering mechanics 120~(2)
  (1994) 368--390.
\newblock \href {https://doi.org/10.1061/(ASCE)0733-9399(1994)120:2(368)}
  {\path{doi:10.1061/(ASCE)0733-9399(1994)120:2(368)}}.

\bibitem{brunton2016discovering}
S.~L. Brunton, J.~L. Proctor, J.~N. Kutz, Discovering governing equations from
  data by sparse identification of nonlinear dynamical systems, Proceedings of
  the national academy of sciences 113~(15) (2016) 3932--3937.
\newblock \href {https://doi.org/10.1073/pnas.1517384113}
  {\path{doi:10.1073/pnas.1517384113}}.

\bibitem{chaturantabut2010nonlinear}
S.~Chaturantabut, D.~C. Sorensen, Nonlinear model reduction via discrete
  empirical interpolation, SIAM Journal on Scientific Computing 32~(5) (2010)
  2737--2764.
\newblock \href {https://doi.org/10.1137/090766498}
  {\path{doi:10.1137/090766498}}.

\bibitem{semaan2017optimal}
R.~Semaan, Optimal sensor placement using machine learning, Computers \& Fluids
  159 (2017) 167--176.
\newblock \href {https://doi.org/10.1016/j.compfluid.2017.10.002}
  {\path{doi:10.1016/j.compfluid.2017.10.002}}.

\bibitem{berkooz1993proper}
G.~Berkooz, P.~Holmes, J.~L. Lumley, The proper orthogonal decomposition in the
  analysis of turbulent flows, Annual review of fluid mechanics 25~(1) (1993)
  539--575.
\newblock \href {https://doi.org/10.1146/annurev.fl.25.010193.002543}
  {\path{doi:10.1146/annurev.fl.25.010193.002543}}.

\bibitem{peherstorfer2018stability}
B.~Peherstorfer, Z.~Drmac, S.~Gugercin, Stability of discrete empirical
  interpolation and gappy proper orthogonal decomposition with randomized and
  deterministic sampling points (2018).
\newblock \href {http://arxiv.org/abs/1808.10473} {\path{arXiv:1808.10473}}.

\bibitem{nakai2020effect}
K.~Nakai, K.~Yamada, T.~Nagata, Y.~Saito, T.~Nonomura, Effect of objective
  function on data-driven sparsesensor optimization (2020).
\newblock \href {http://arxiv.org/abs/2007.05377} {\path{arXiv:2007.05377}}.

\bibitem{ye2018complexity}
L.~Ye, S.~Roy, S.~Sundaram, On the complexity and approximability of optimal
  sensor selection for kalman filtering, in: 2018 Annual American Control
  Conference (ACC), IEEE, 2018, pp. 5049--5054.
\newblock \href {https://doi.org/10.23919/ACC.2018.8431016}
  {\path{doi:10.23919/ACC.2018.8431016}}.

\bibitem{joshi2009sensor}
S.~Joshi, S.~Boyd, Sensor selection via convex optimization, IEEE Transactions
  on Signal Processing 57~(2) (2009) 451--462.
\newblock \href {https://doi.org/10.1109/TSP.2008.2007095}
  {\path{doi:10.1109/TSP.2008.2007095}}.

\bibitem{nonomura2020randomized}
T.~Nonomura, S.~Ono, K.~Nakai, Y.~Saito, Randomized subspace newton convex
  method applied to data-driven sensor selection problem (2020).
\newblock \href {http://arxiv.org/abs/2009.09315} {\path{arXiv:2009.09315}}.

\bibitem{manohar2018data}
K.~Manohar, B.~W. Brunton, J.~N. Kutz, S.~L. Brunton, Data-driven sparse sensor
  placement for reconstruction: Demonstrating the benefits of exploiting known
  patterns, IEEE Control Systems Magazine 38~(3) (2018) 63--86.
\newblock \href {https://doi.org/10.1109/MCS.2018.2810460}
  {\path{doi:10.1109/MCS.2018.2810460}}.

\bibitem{saito2019determinantbased}
Y.~Saito, T.~Nonomura, K.~Yamada, K.~Asai, Y.~Sasaki, D.~Tsubakino,
  Determinant-based fast greedy sensor selection algorithm (2019).
\newblock \href {http://arxiv.org/abs/1911.08757} {\path{arXiv:1911.08757}}.

\bibitem{dhingra2014admm}
N.~K. Dhingra, M.~R. Jovanovi{\'c}, Z.-Q. Luo, An admm algorithm for optimal
  sensor and actuator selection, in: 53rd IEEE Conference on Decision and
  Control, IEEE, 2014, pp. 4039--4044.
\newblock \href {https://doi.org/10.1109/CDC.2014.7040017}
  {\path{doi:10.1109/CDC.2014.7040017}}.

\bibitem{nagata2020data-driven}
T.~Nagata, T.~Nonomura, K.~Nakai, K.~Yamada, Y.~Saito, S.~Ono, Data-driven
  sparse sensor placement based on a-optimal design of experiment with admm
  (2020).
\newblock \href {http://arxiv.org/abs/2010.09329} {\path{arXiv:2010.09329}}.

\bibitem{o2016distributed}
M.~O'Connor, W.~B. Kleijn, T.~Abhayapala, Distributed sparse mvdr beamforming
  using the bi-alternating direction method of multipliers, in: 2016 IEEE
  International Conference on Acoustics, Speech and Signal Processing (ICASSP),
  IEEE, 2016, pp. 106--110.
\newblock \href {https://doi.org/10.1109/ICASSP.2016.7471646}
  {\path{doi:10.1109/ICASSP.2016.7471646}}.

\bibitem{castro2013robustness}
R.~Castro-Triguero, S.~Murugan, R.~Gallego, M.~I. Friswell, Robustness of
  optimal sensor placement under parametric uncertainty, Mechanical Systems and
  Signal Processing 41~(1-2) (2013) 268--287.
\newblock \href {https://doi.org/10.1016/j.ymssp.2013.06.022}
  {\path{doi:10.1016/j.ymssp.2013.06.022}}.

\bibitem{liu2016sensor}
S.~Liu, S.~P. Chepuri, M.~Fardad, E.~Ma{\c{s}}azade, G.~Leus, P.~K. Varshney,
  Sensor selection for estimation with correlated measurement noise, IEEE
  Transactions on Signal Processing 64~(13) (2016) 3509--3522.
\newblock \href {https://doi.org/10.1109/TSP.2016.2550005}
  {\path{doi:10.1109/TSP.2016.2550005}}.

\bibitem{ucinski2020d}
D.~Uci{\'n}ski, D-optimal sensor selection in the presence of correlated
  measurement noise, Measurement 164 (2020) 107873.
\newblock \href {https://doi.org/10.1016/j.measurement.2020.107873}
  {\path{doi:10.1016/j.measurement.2020.107873}}.

\bibitem{nankai2019linear}
K.~Nankai, Y.~Ozawa, T.~Nonomura, K.~Asai, Linear reduced-order model based on
  piv data of flow field around airfoil, TRANSACTIONS OF THE JAPAN SOCIETY FOR
  AERONAUTICAL AND SPACE SCIENCES 62~(4) (2019) 227--235.
\newblock \href {https://doi.org/10.2322/tjsass.62.227}
  {\path{doi:10.2322/tjsass.62.227}}.

\bibitem{saito2020data}
Y.~Saito, T.~Nonomura, K.~Nankai, K.~Yamada, K.~Asai, Y.~Sasaki, D.~Tsubakino,
  Data-driven vector-measurement-sensor selection based on greedy algorithm,
  IEEE Sensors Letters 4 (2020).
\newblock \href {https://doi.org/10.1109/LSENS.2020.2999186}
  {\path{doi:10.1109/LSENS.2020.2999186}}.

\bibitem{noaa}
NOAA/OAR/ESRL,
  \href{https://www.esrl.noaa.gov/psd/data/gridded/data.noaa.oisst.v2.html}{Noaa
  optimal interpolation (oi) sea surface temperature (sst) v2} (July 2019).
\newline\urlprefix\url{https://www.esrl.noaa.gov/psd/data/gridded/data.noaa.oisst.v2.html}

\bibitem{summers2015submodularity}
T.~H. Summers, F.~L. Cortesi, J.~Lygeros, On submodularity and controllability
  in complex dynamical networks, IEEE Transactions on Control of Network
  Systems 3~(1) (2015) 91--101.
\newblock \href {https://doi.org/10.1109/TCNS.2015.2453711}
  {\path{doi:10.1109/TCNS.2015.2453711}}.

\end{thebibliography}

\end{document}